\documentclass[11pt]{amsart}
\usepackage{amsmath,amssymb,amscd}
\usepackage{mathrsfs}
\usepackage{graphicx}
\usepackage[all]{xy}

\title{Projective models of the twistor spaces of Joyce metrics}

\author{Nobuhiro Honda$^{\dag}$}
\thanks
{$^{\dag}$Partially supported by
Research Fellowships of the 
Japan Society for the Promotion
of Science for Young Scientists.
}
\date{}

\newcommand{\ol}{\overline}

\newcommand{\lra}{\longrightarrow}

\newcommand{\set}{\,|\,}
\newcommand{\proofend}{\hfill$\square$}

\newtheorem{prop}{Proposition}[section]
\newtheorem{lemma}[prop]{Lemma}

\newtheorem{thm}[prop]{Theorem}

\newtheorem{definition}[prop]{Definition}

\begin{document}

\begin{abstract}
We provide a simple algebraic construction of the twistor spaces of arbitrary Joyce's self-dual metrics on the 4-manifold $\mathscr H^2\times T^2$ that extend  smoothly to $n\mathbf{CP}^2$, the connected sum of complex projective planes. 
Indeed, we explicitly realize projective models of the twistor spaces of arbitrary Joyce metrics on $n\mathbf{CP}^2$ in a $\mathbf{CP}^4$-bundle over $\mathbf{CP}^1$, and show that they contain the twistor spaces of $\mathscr H^2\times T^2$ as dense non-Zariski open subsets.
In particular, we see that the last non-compact twistor spaces can be realized in  rank-4 vector bundles over $\mathbf{CP}^1$ by quite simple defining equations.
\end{abstract}

\maketitle
\section{Introduction}
In the paper \cite{J95}, D.\,Joyce found a beautiful series of self-dual conformal classes on the connected sum of complex projective planes that are invariant under effective $T^2$-actions.
His construction consists of two parts:

\begin{itemize}
\item
%
 First, explicit self-dual metrics on the non-compact 4-manifold $\mathscr H^2\times T^2$ are written down in coordinate. These  metrics are invariant under the natural $T^2$-action, and depend on $(n+2)$ real numbers and also $(n+2)$ vectors in $\mathbf R^2$.

\item
 Then it is shown that if the latter vectors take special integer values determined by an effective $T^2$-action on $n\mathbf{CP}^2$, then the metrics can be extended smoothly to a compactification that is diffeomorphic to $n\mathbf{CP}^2$.
\end{itemize}
In the paper he asks a simple construction of the twistor spaces \cite[p.\,548]{J95}.
The purpose of this paper is to give a simple algebraic presentation of the twistor spaces of these Joyce metrics on $\mathscr H^2\times T^2$ that extend to $n\mathbf{CP}^2$.
In the course of the construction we also provide projective models of the twistor spaces for arbitrary Joyce metrics on $n\mathbf{CP}^2$.

A main tool of the present investigation is our recent result on linear systems on the twistor spaces of Joyce metrics on $n\mathbf{CP}^2$ \cite{Hon08-1}.
There, for arbitrary $U(1)$-subgroup of $T^2$ that fixes a $T^2$-invariant sphere in $n\mathbf{CP}^2$, we have concretely found a linear system on the twistor space whose meromorphic map can be regarded as a quotient map with respect to the natural $\mathbf C^*$-action corresponding to the $U(1)$-action.
In this paper, by taking two different such $U(1)$-subgroups in $T^2$ and adjusting one of  the resulting two linear systems slightly, we obtain a single, $T^2$-invariant linear system.
It is possible to give generators of the last linear system and their relations in concrete forms.
Hence we can obtain  defining equations of the image of the associated meromorphic map in explicit forms.
We show that by an appropriate (and explicit) blowing-up, the images are transformed into 3-folds 
embedded in  $\mathbf{CP}^4$-bundles over $\mathbf{CP}^1$. 
The defining equations and geometric structure of the images take most simple form in this situation:
the defining equations become two quadratic polynomials defined over the polynomial ring $\mathbf R[\lambda]$, where the polynomials are uniquely determined by the invariants involved in the construction of the Joyce metrics, and the 3-folds are fibered over $\mathbf{CP}^1$ by compact toric surfaces.

Although the last 3-folds in  $\mathbf{CP}^4$-bundles over $\mathbf{CP}^1$ are meromorphic images of the original twistor spaces of Joyce metrics on $n\mathbf{CP}^2$, the meromorphic maps can become $d:1$ with $d>1$; namely they are not necessarily bimeromorphic.
We give a simple criterion for the map to be bimeromorphic, and also remark that we can always find the cases for the criterion to be satisfied.
This way we obtain projective models of arbitrary Joyce metrics on $n\mathbf{CP}^2$.

Furthermore, by making use of the explicitness of our meromorphic maps,   we show that the bimeromorphic maps from the twistor spaces to the projective models never contract any two-dimensional orbits of $(\mathbf C^*)^2$-action on the twistor spaces, where the $(\mathbf C^*)^2$-action is the complexification of the $T^2$-action.
This in particular means that the meromorphic map is biholomorphic on the open subset which is the twistor space of $\mathscr H^2\times T^2$.
The explicitness also enable us to conclude that the last twistor spaces are holomorphically embedded in  rank-4 complex vector bundles over $\mathbf{CP}^1$.
Thus unlike the twistor spaces on the whole of $n\mathbf{CP}^2$, the twistor spaces on $\mathscr H^2\times T^2$ can always be embedded  in projective spaces. 

It should be mentioned relationship between our previous paper  \cite{Hon06-1}.
In the paper we have constructed not only projective models of the twistor spaces of Joyce metrics on $4\mathbf{CP}^2$ but also their biholomorphic models. 
As is evident from the construction there, to obtain the biholomorphic models from the projective models requires quite complicated and partly ad-hoc operations.
Hence, the present result is stronger in that it concerns {\em arbitrary} Joyce metrics, but weaker in that it does not give the biholomorphic models of the twistor spaces {\em on the whole of}  $n\mathbf{CP}^2$.
We remark that, at the metric level also, the situation looks similar in that  the self-dual metrics are explicitly written down on $\mathscr H^2\times T^2$ but not  on the whole of $n\mathbf{CP}^2$.

{\bf Acknowledgement.} I would like to thank Junjiro\,Noguchi whose question  made me notice that the operations  from projective models to the biholomorphic models given in \cite{Hon06-1} do not touch 2-dimensional orbits at all.

{\bf Notations.}
For a twistor space $Z$ of a self-dual 4-manifold $M$, $F$ denotes the canonical square root of the anticanonical line bundle of $Z$.
If a group $G$ is acting on $M$ preserving the self-dual structure, $G$ acts naturally on $Z$ holomorphically.
The last $G$-action also lifts on the space $H^0(Z,mF)$ of sections of $mF=F^{\otimes m}$.
The subspace consisting of all $G$-invariant sections is denoted by $H^0(Z,mF)^G$.
If $V$ is a linear subspace of $H^0(Z, mF)$, $|V|$ denotes the linear subsystem of $|mF|$ whose members are defined by $s\in V$.
If $V=H^0(Z,mF)^G$, then $|V|$ is denoted by $|mF|^G$.
For $0\neq s\in V$, $(s)$ denotes the divisor defined as the zero-locus of $s$.
If $W$ is a $\mathbf C$-vector space, $\mathbf P^{\vee}W$ denotes the dual complex projective space.

\section{Algebraic construction of the twistor spaces}
First, for the purpose of explaining the results of Joyce and Fujiki in \cite{J95, F00} and fixing some notations, we briefly recall basic properties of $T^2$-actions on $n\mathbf{CP}^2$.
Fix an effective $T^2$-action on $n\mathbf{CP}^2$.
Then the set of points on $n\mathbf{CP}^2$ having non-trivial isotropy subgroup consists of  $k:=n+2$ smooth $T^2$-invariant 2-spheres $B_i$, $1\le i\le k$.
The union of these spheres form a cycle, in the sense that $B_i$ and $B_{i+1}$ intersect transversally at a single point $p_i$, where we interpret $B_{k+1}=B_1$, and that other intersections are empty.
$T^2$ acts freely  on the complement  of this cycle, and there exists a $T^2$-equivariant diffeomorphism  
\begin{align}
n\mathbf{CP}^2\backslash (B_1\cup\cdots \cup B_{k})\simeq \mathscr H^2\times T^2
\end{align}
where $T^2$ is acting on $\mathscr H^2\times T^2$ in an obvious way.
Thus $\mathscr H^2\times T^2$ is contained in $n\mathbf{CP}^2$ as a $T^2$-invariant dense open subset.
On each $B_i$,  a subgroup of form 
\begin{align}
G(n_i,n'_i):=\{(s,t)\in T^2=U(1)\times U(1)\set s^{n_i}t^{n'_i}=1\}
\end{align}
is acting trivially, where the primitive vectors $(n_i,n'_i)\in\mathbf Z^2$ are uniquely determined up to simultaneous inversion of the sign.
By the standard $SL(2,\mathbf Z)$-action, one can always normalize the set $\{(n_i,n'_i)\set 1\le i\le k\}$ so as to satisfy the following conditions \cite[Prop.\,3.1.1]{J95}: 
\begin{align}
&(n_1,n'_1)=(0,1),\,\,(n_k,n'_k)=(1,0),\label{normal1}\\
&n_i>0 {\text{ for }} i>0,\label{normal2}\\
&n_in'_{i+1}-n_{i+1}n'_i=-1  {\text{ for }} 1\le i<k.\label{normal3}
\end{align}
Conversely, any set $\{(n_i,n'_i)\in\mathbf Z^2\set 1\le i\le k\}$  satisfying \eqref{normal1}--\eqref{normal3} determines an effective $T^2$-action on $n\mathbf{CP}^2$ \cite[Prop.\,3.1.1]{J95}.

In \cite{J95} Joyce constructed a series of $T^2$-invariant self-dual metrics on $\mathscr H^2\times T^2$ that can be written down in coordinates. 
His metrics involve various parameters:
 an integer $n\ge 0$,
different real numbers $\lambda_1,\lambda_2,\cdots,\lambda_{k}\in\mathbf R\cup\{\infty\}$, and 
 vectors $v_1,v_2,\cdots,v_{k}\in\mathbf R^2$.
 He further showed that if  $v_i=(n_i,n'_i)$ for some set $\{(n_i,n'_i)\set 1\le i\le k\}$ satisfying  \eqref{normal1}--\eqref{normal3}, then the metric on $\mathscr H^2\times T^2$ extends smoothly to $n\mathbf{CP}^2$ in such a way that the $T^2$-action on $n\mathbf{CP}^2$  is the one determined by $\{(n_i,n'_i)\set 1\le i\le k\}$ (\cite[Theorem 3.3.1]{J95}).

In the following we fix any $n\ge 0$, $\lambda_1,\cdots,\lambda_k\in\mathbf R\cup\{\infty\}$ and  $\{(n_i,n'_i)\set 1\le i\le k\}$ satisfying  \eqref{normal1}--\eqref{normal3}, and let $g$ be the corresponding Joyce metric on $n\mathbf{CP}^2$.
The restriction of $g$ to the open dense  subset $\mathscr H^2\times T^2$  is denoted by $g^{\circ}$.
Let $Z$ and $Z^{\circ}$ be the twistor spaces of $g$ and $g^{\circ}$ respectively.
$Z^{\circ}$ is a non-Zariski open dense subset of $Z$.
$Z$ is equipped with a holomorphic $(\mathbf C^*)^2$-action obtained from  the given $T^2$-action on $n\mathbf{CP}^2$, and $Z^{\circ}$ is $(\mathbf C^*)^2$-invariant in $Z$.

The structure of $Z$ is extensively investigated by A.\,Fujiki in \cite{F00}, by effectively making use of the $(\mathbf C^*)^2$-action.
We collect some important properties among them we will need.

\begin{prop}\label{prop-Fujiki1}
(i) The linear system $|F|^{T^2}$ is a pencil whose general members are mutually isomorphic non-singular toric surfaces.
(ii) The set of rays of a fan corresponding to the toric surface in (i) is given by $\{\mathbf R_{+}v_i,\mathbf R_{+}(-v_i)\set 1\le i\le k\}$.
(iii) Letting $C_i$ and $\ol C_i$ be the closures of 1-dimensional orbits corresponding to  $\mathbf R_{+}(v_i)$ and $\mathbf R_{+}(-v_i)$ respectively, so that $C:=\sum_{1\le i\le k}(C_i+\ol C_i)$ is a cycle of anticanonical curve of the toric surface, $C$ coincides with the base locus of $|F|^{T^2}$.
(iv) The real structure of $Z$ exchanges $C_i$ and $\ol C_i$.
(v) The singular members of $|F|^{T^2}$ consist of $k$ divisors of the form $S_i^++S_i^-$, $1\le i\le k$, where $S_i^+$ and $S_i^-$ are mutually conjugate non-singular toric surfaces.
\end{prop}

If $Z$ is not a LeBrun twistor space constructed in \cite{LB91}, we can replace $|F|^{T^2}$ by just $|F|$.
He further showed  that after a possible renumbering the intersections $L_i:=S_i^+\cap S_i^-$ are precisely the twistor lines over the $T^2$-fixed points $p_i=B_i\cap B_{i+1}$.
Thus $S_i^+$ and $S_i^-$ divide the cycle $C$ into connected `halves'.
By $T^2$-equivariancy, the restriction of the twistor fibration $Z\to n\mathbf{CP}^2$ to $C$ gives a unramified double covering map  onto $B_1\cup\cdots\cup B_k$, where $C_i$ and $\ol C_i$ are mapped to $B_i$.
By Prop.\,\ref{prop-Fujiki1} the union of all 0- and  1-dimensional orbits of the $(\mathbf C^*)^2$-action on $Z$ is given by
$C\cup(L_1\cup\cdots\cup L_k).$
We set its complement as
\begin{align}\label{U}
U:=Z\backslash (C\cup(L_1\cup\cdots\cup L_k)).
\end{align}
Then $U$ is of course a Zariski open subset of $Z$, acted freely by $(\mathbf C^*)^2$.
Since $Z^{\circ}$ is the complement in $Z$ of all twistor lines over $B_1\cup\cdots\cup B_k$, we have $Z^{\circ}\subset U$.
The complement $U\backslash Z^{\circ}$ is an open dense subset of the union of all twistor lines over $(B_1\cup\cdots\cup B_k)$.
Thus $U$ is slightly larger than the twistor space $Z^{\circ}$ of $\mathscr H^2\times T^2$.

We also recall a result in \cite{Hon08-1}.
For this, as in \cite{Hon08-1}, for a positive integer $m$, let $V_m\subset H^0(Z,mF)$ be the subspace generated by the image of the natural map $(H^0(Z,F)^{T^2})^{\times m}\to H^0(Z,mF)$ given by $(s_1,\cdots,s_m)\mapsto s_1\cdot\cdots\cdot s_m$.
Since $V_1=H^0(Z,F)^{T^2}\simeq\mathbf C^2$, we have $\dim V_m=m+1$.
Further, as showed in \cite[Prop.\,2.11]{Hon08-1}, $V_m=H^0(Z,mF)^{T^2}$ holds.
Let $\Lambda_m\subset\mathbf P^{\vee}V_m$ be the rational normal curve determined via the natural map $V_1\to V_m$ defined as $s\mapsto s^m$.
When $m$ varies, all these rational normal curves are naturally identified via the isomorphisms $\Lambda_1=\mathbf P^{\vee}V_1\simeq\Lambda_m$.
Also, we put $G_i=G(n_i,n'_i)$.
Then we have the following (\cite[Prop.\,2.5]{Hon08-1}).

\begin{prop}\label{prop-mt1}
For any $1\le i\le k$, we can explicitly find an integer $m_i$ and a divisor $Y_i\in |m_iF|$ satisfying the following.
(i) The linear system generated by  $|V_{m_i}|$, $Y_i$ and $\ol Y_i$ coincides with $|m_iF|^{G_i}$.
(ii) $\dim |m_iF|^{G_i}=m_i+2$.
(iii) All irreducible components of $Y_i$ are of degree one,
(iv) $Y_i$ does not contain $S_j^+$ and $S_j^-$ at the same time for any $1\le j\le k$.
\end{prop}
In \cite{Hon08-1},  $m_i$ and  $Y_i$ are denoted by just $m$ and $Y$ respectively, since we were fixing one of the $T^2$-invariant spheres.
As explained in \cite{Hon08-1}, once a $T^2$-action on $n\mathbf{CP}^2$ and $1\le i\le k$ are given, it is  easy to determine $m_i$ and $Y_i$. Indeed the computations have been given in a simple algorithmic form.
For later reference, for each $1\le i\le k$, we also prepare the following notations.

\begin{definition}
{\em 
(i) We put $W_{m_i}:=H^0(Z,m_iF)^{G_i}.$
(ii) Take a (non-real) section $y_i\in W_{m_i}$ such that   $(y_i)=Y_i$.
The conjugate section  is denoted by $\ol y_i\,(\in |W_{m_i}|)$, which satisfies $(\ol y_i)=\ol Y_i$.
(iii) 
Take a real section $u_i\in H^0(Z,F)$ such that  $(u_i)=S_i^++S_i^-$.
}
\end{definition}
\noindent
By Prop.\,\ref{prop-mt1}, as  basis of $V_{m_i}$ and $W_{m_i}$ we can explicitly take 
\begin{equation}\label{base1}
V_{m_i}=\langle u_1^{m_i}, u_1^{m_i-1}u_2,\cdots,u_2^{m_i}\rangle
\subset W_{m_i}=V_{m_i}\oplus \langle y_i,\ol y_i\rangle.
\end{equation}
Here, we are choosing $\{u_1,u_2\}$ as a basis of $H^0(Z, F)^{T^2}$.
By \cite[Theorem 9.1]{F00}, we can write
\begin{align}\label{lcb}
u_{\alpha}=u_1-\lambda_{\alpha}u_2,\hspace{2mm}{\rm{for}}\hspace{2mm}3\le\alpha\le k,
\end{align}
where the invariant $\{\lambda_1,\cdots,\lambda_k\}$ is normalized so as to satisfy $\lambda_1=\infty,\lambda_2=0$ and either $0<\lambda_2<\lambda_3<\cdots<\lambda_k$ or $0>\lambda_2>\lambda_3>\cdots>\lambda_k$ holds.

Our goal is to realize the manifold $U$ in \eqref{U} explicitly  in certain $\mathbf{CP}^4$-bundle over $\mathbf{CP}^1$.
For this, we first choose arbitrary two components $B_i$ and $B_j$ $(1\le i<j\le k)$ among the $T^2$-invariant spheres in $n\mathbf{CP}^2$. 
We may suppose $m_i\ge m_j$.
(If not, it is enough to replace the role of $i$ and $j$ in the following construction.)
For simplicity of notation, we put
\begin{align}
\mu:=m_i-m_j\,(\ge 0)
\end{align}
Let $V_{\mu}\cdot W_{m_j}$ be the image of the canonical map from  $V_{\mu}\otimes W_{m_j}$ (the tensor product of vector spaces in the usual sense) to $ H^0(Z,m_iF)$.
(In the case of $m_i=m_j$, $V_0$ is interpreted as  $\mathbf C$, so that $V_{\mu}\cdot W_{m_j}=W_{m_j}$.)
 $V_{\mu}\cdot W_{m_j}$ is a $T^2$-invariant subspace, acted trivially by $G_j$. Then we consider the vector subspace
\begin{align}\label{w1}
W_{m_i}+ V_{\mu}\cdot W_{m_j}\subset H^0(Z,m_iF),
\end{align}
which is also $T^2$-invariant.
The dimension of  \eqref{w1} can be readily computed by the following 

\begin{prop}
$W_{m_i}\cap (V_{\mu}\cdot W_{m_j})=V_{m_i}.$

\end{prop}

\noindent 
Proof.
For an inclusion `$\supset$', since $W_{m_i}\supset V_{m_i}$, it suffices to show $V_{\mu}\cdot W_{m_j}\supset V_{m_i}$.
But this is obvious since $W_{m_j}\supset V_{m_j}$ and $V_{\mu}\cdot V_{m_j}=V_{m_i}$. Hence we obtain the inclusion.
Another inclusion is also immediate by $W_{m_i}=H^0(Z,m_iF)^{G_i}$, $V_{\mu}\cdot W_{m_j}\subset H^0(Z, m_iF)^{G_j}$, and  $H^0(Z,m_iF)^{G_i}\cap H^0(Z,m_iF)^{G_j}=H^0(Z,m_iF)^{T^2}=V_{m_i}$.
Thus we obtain the claim of the proposition.
\proofend

\vspace{2mm}
By this proposition, for a natural basis of $W_{m_i}+ V_{\mu}\cdot W_{m_j}$, it is enough to add the following $2(\mu+1)$ sections
\begin{align}\label{base2}
u_1^{\mu_1}u_2^{\mu_2} y_j,\,\,\,u_1^{\mu_1}u_2^{\mu_2}\ol y_j\
\hspace{2mm}{\text{while}}\hspace{2mm}\mu_1+\mu_2=\mu
\end{align}
to the basis \eqref{base1} of $W_{m_i}$.
From this we obtain
$\dim (W_{m_i}+ V_{\mu}\cdot W_{m_j})=3m_i-2m_j+5.$
(This value itself is not so important.)
Let $(\zeta_1,\zeta_2,\cdots,\zeta_{N})$, $N=5+3m_i-2m_j$, be a homogeneous coordinate on $ \mathbf P^{\vee}\left(W_{m_i}+ V_{\mu}\cdot W_{m_j}\right)$.
Let 
\begin{align}
\Phi:Z\to \mathbf P^{\vee}\left(W_{m_i}+ V_{\mu}\cdot W_{m_j}\right)
\end{align}
be the meromorphic map associated to the linear system $|W_{m_i}+ V_{\mu}\cdot W_{m_j}|$. 
(So $\Phi$ depends on the choice of $1\le i<j\le k$.)
$\Phi$ is a $T^2$-equivariant map. 
By the choice of the basis of $W_{m_i}+ V_{\mu}\cdot W_{m_j}$ in \eqref{base1} and \eqref{base2}, we can suppose that  $\Phi$ is explicitly given by 
\begin{align}\label{hc1}
\zeta_1=y_i,\,\zeta_2=\ol y_i, 
\end{align}
\begin{align}\label{hc2}
\zeta_3=u_1^{\mu}y_j,\,\zeta_4=u_1^{\mu}\ol y_j,
\zeta_5=u_1^{\mu-1}u_2y_j,\,\zeta_4=u_1^{\mu-1}u_2\ol y_j, \cdots,
\zeta_{2\mu+3}=u_2^{\mu}y_j,\,\zeta_{2\mu+4}=u_2^{\mu}\ol y_j
\end{align}
and
\begin{align}\label{hc3}
\zeta_{N-m_i}=u_1^{m_i},\, \zeta_{N-m_i+1}=u_1^{m_i-1}u_2,\cdots,\,\zeta_{N}=u_2^{m_i}.
\end{align}
(Note $N-m_i=2\mu+5$.)
Then the rational normal curve $\Lambda_{m_i}\subset \mathbf P^{\vee}V_{m_i}$ is explicitly given by
\begin{align}
\Lambda_{m_i}=\{(\zeta_{N-m_i},\zeta_{N-m_i+1},\cdots,\zeta_N)=(1,\lambda,\lambda^2,\cdots,\lambda^{m_i})\set\lambda\in\mathbf C\}\cup
\{(0,0,\cdots,0,1)\}.
\end{align}
On the other hand we have a basic commutative diagram of meromorphic maps
\begin{equation}\label{cd1}
 \CD
Z@>{\Phi}>>\mathbf P^{\vee}\left(W_{m_i}+ V_{\mu}\cdot W_{m_j}\right)\\
 @V\Psi VV @VV{p}V\\
\mathbf P^{\vee}V_1@>>>\mathbf P^{\vee}V_{m_i},\\
 \endCD
 \end{equation}
where $\Psi$ is the meromorphic map associated to the pencil $|V_1|=|F|^{T^2}$, $p$ is the projection induced by the inclusion $V_{m_i}\subset W_{m_i}$, and $\mathbf P^{\vee}V_1\to\mathbf P^{\vee}V_{m_i}$ is an inclusion onto the rational normal curve $\Lambda_{m_i}$ as before.
This diagram is clearly $T^2$-equivariant, where $T^2$ is acting trivially on the bottom two projective spaces.
In the above homogeneous coordinates,
$p$ is written as $(\zeta_{1},\cdots,\zeta_N)\mapsto(\zeta_{N-m_i},\cdots,\zeta_N)$.
If we blow-up the indeterminacy locus of   $p$ explicitly given by
\begin{align}\label{center1}
\mathbf P_{\infty}:=\{\zeta_{N-m_i}=\cdots=\zeta_N =0\},
\end{align}
then we obtain the total space of the $\mathbf P^{2(\mu+2)}$-bundle
\begin{align}\label{bu1}
\mathbf P(\mathscr O(1)^{\oplus 2(\mu+2)}\oplus\mathscr O)\lra \mathbf P^{\vee}V_{m_i}.
\end{align}
Since $\mathbf P_{\infty}$ is a $T^2$-invariant subspace, the bundle \eqref{bu1} still possesses a natural $T^2$-action, and all fibers are $T^2$-invariant.
Let $\tilde \Phi:Z\to \mathbf P(\mathscr O(1)^{\oplus 2(\mu+2)}\oplus\mathscr O)$ be the meromorphic map obtained as the composition of $\Phi$ and the inverse of the blowing-up at $\mathbf P_{\infty}$.
By the diagram \eqref{cd1}, the image $\tilde\Phi(Z)$ is contained in the restriction of the bundle \eqref{bu1} onto the curve $\Lambda_{m_i}$, and the latter is of course isomorphic to the bundle 
\begin{align}\label{bd2}
\mathbf P(\mathscr O(m_i)^{\oplus 2(\mu+2)}\oplus\mathscr O)\lra \Lambda_{m_i}\simeq\mathbf{CP}^1.
\end{align}
$T^2$ is still acting on this bundle, keeping the fibers invariant, and $\tilde\Phi$ is also $T^2$-equivariant.
We are going to determine the image $\tilde\Phi(Z)$ in this bundle concretely.
For this, we put
\begin{align}\label{fc1}
\xi_{\alpha}=\frac{\zeta_{\alpha}}{\zeta_{N-m_i}}, \,\,1\le \alpha\le 2(\mu+2)
\end{align}
and use $(\xi_1,\cdots,\xi_{2(\mu+2)})$ as a non-homogeneous fiber coordinate on $\mathbf P(\mathscr O(1)^{\oplus 2(\mu+2)}\oplus\mathscr O)\lra \mathbf P^{\vee}V_{m_i}$
valid on the Zariski open set $\{\zeta_{N-m_i}\neq0\}\subset\mathbf P^{\vee}V_{m_i}$.
Namely, $\xi_{\alpha}$ is a fiber coordinate on the $\alpha$-th line bundle $\mathscr O(1)$ (or $\mathscr O(m_i)$ over $\Lambda_{m_i}$) of the bundle $\mathscr O(1)^{\oplus2(\mu+1)}$ (or $\mathscr O(m_i)^{\oplus2(\mu+1)}\to\Lambda_{m_i}$). The last vector bundle is canonically embedded in the bundle \eqref{bu1} (or \eqref{bd2}) as a $T^2$-invariant Zariski open subset.
We also note that  $\xi_{1}$ and $\xi_2$ are acted trivially by $G_i$, while 
$\xi_3,\xi_4,\cdots,\xi_{2\mu+4}$ are acted trivially by $G_j$.
Moreover the product $\xi_{2\alpha-1}\xi_{2\alpha}$ is $T^2$-invariant for any $1\le \alpha\le \mu+2$.
Furthermore, let
\begin{align}
Y_{i}=\sum_{1\le \beta\le k}(l_{i\beta}^+S_{\beta}^++l_{i\beta}^-S_{\beta}^-)\hspace{2mm}{\text{and}}\hspace{2mm}
Y_{j}=\sum_{1\le \beta\le k}(l_{j\beta}^+S_{\beta}^++l_{j\beta}^-S_{\beta}^-)
\end{align}
be the decomposition into irreducible components, where $l_{i\beta}^+l_{i\beta}^-=l_{j\beta}^+l_{j\beta}^-=0$ hold by Prop.\,\ref{prop-mt1} (iv).
Also we define non-negative integers $l_{i\beta}$ and $l_{j\beta}$ by $l_{i\beta}=l_{i\beta}^++l_{i\beta}^-$ and $l_{j\beta}=l_{j\beta}^++l_{j\beta}^-$.
(These numbers only depend on a $T^2$-action on $n\mathbf{CP}^2$ and the choice of $i$ and $j$.
As is already noted, all these are explicitly and readily computable numbers.)

\begin{prop}\label{prop-image1}
In the above coordinate, the image  $\tilde\Phi(Z)$ is contained in an algebraic variety defined by the following equations:
\begin{align}\label{mt1}
\xi_1\xi_2=c_1\prod_{2\le \beta\le k}(\lambda-\lambda_{\beta})^{l_{i\beta}},
\end{align}
\begin{align}\label{mt2}
\xi_{2\alpha-1}\xi_{2\alpha}=c_{\alpha}\lambda^{2(\alpha-2)}\prod_{2\le \beta\le k}(\lambda-\lambda_{\beta})^{l_{j\beta}},\hspace{2mm}2\le\alpha\le \mu+2,
\end{align}
where $c_1$ and $c_{\alpha}$ are non-zero real constants.
\end{prop}

\noindent Proof.
For $\gamma=i$ and $j$ we have 
\begin{align}
Y_{\gamma}+\ol Y_{\gamma}=\sum_{1\le \beta\le k}l_{\gamma\beta}(S_{\beta}^++S_{\beta}^-).
\end{align}
Therefore, since $(u_{\beta})=S_{\beta}^++S_{\beta}^-$, we obtain, for sections
\begin{align}\label{dcp3}
y_{\gamma}\ol y_{\gamma}&=c_{\gamma}u_1^{l_{\gamma1}}u_2^{l_{\gamma2}}\cdots u_k^{l_{\gamma k}}
\end{align}
for some non-zero real constant $c_{\gamma}$.
By \eqref{lcb}, the right-hand side can be written as 
\begin{align}
c_{\gamma}u_1^{l_{\gamma1}}u_2^{l_{\gamma2}}(u_1-\lambda_3u_2)^{l_{\gamma3}}\cdots(u_1-\lambda_ku_2)^{l_{\gamma k}}
\end{align}
Then by \eqref{hc1} and \eqref{hc3}, we have
 \begin{align}
 \xi_1\xi_2=\frac{y_i\ol y_i}{u_1^{2m_i}}=\frac{c_{\gamma}u_1^{l_{\gamma1}}u_2^{l_{\gamma2}}(u_1-\lambda_3u_2)^{l_{\gamma3}}\cdots(u_1-\lambda_ku_2)^{l_{\gamma k}}}{u_1^{2m_i}}.
 \end{align}
Hence by putting $\lambda=u_2/u_1$ and recalling $\sum_{1\le \beta\le k}l_{i\beta}=2m_{i}$, we obtain the equation \eqref{mt1}.
Other equations \eqref{mt2} can be obtained analogously by using \eqref{hc1} and \eqref{hc2}.
\proofend

\vspace{2mm}
The right-hand side of \eqref{mt1} and \eqref{mt2} are of degree $2m_i-l_{i1}$ and $2m_j-l_{j1}+2(\alpha-2)$ respectively.
(If we use a `general' basis of $H^0(F)^{T^2}$ instead of $\{u_1,u_2\}$, the degree of these polynomials become $2m_i$ and $2m_j+2(\alpha-2)$ respectively.)
Counting the number of the equations,  the subvariety defined by \eqref{mt1} and \eqref{mt2} is $(\mu+3)$-dimensional.
Hence if $\mu>0$ the subvariety cannot be the image of $\tilde\Phi$.
This can be remedied by the following

\begin{prop}
There exists a rank-4 subbundle $\mathscr E$ of  $\mathscr O(m_i)^{\oplus 2(\mu+2)}\to\Lambda_{m_i}$
 satisfying the following properties.
 (i) The image $\tilde\Phi(Z)$ is contained in the subbundle $\mathbf P(\mathscr E\oplus\mathscr O)$ of the bundle \eqref{bd2}.
(ii) $\mathscr E$ is isomorphic to $\mathscr O(m_i)^{\oplus 2}\oplus\mathscr O(m_j)^{\oplus 2}$.
(iii) For a fiber coordinate $(\xi_1,\xi_2)$ on $\mathscr O(m_i)^{\oplus 2}$ and $(\xi_3,\xi_4)$ on $\mathscr O(m_j)^{\oplus 2}$, $\tilde{\Phi}(Z)$ is contained in a subvariety defined by
\begin{align}\label{de2}
\xi_1\xi_2=c_1\prod_{2\le \beta\le k}(\lambda-\lambda_{\beta})^{l_{i\beta}},
\hspace{3mm}
\xi_{3}\xi_{4}=c_{2}\prod_{2\le \beta\le k}(\lambda-\lambda_{\beta})^{l_{j\beta}}
\end{align}
for non-zero constants $c_1$ and $c_2$.
\end{prop}

\noindent
Proof.
Let $(\xi_1,\xi_2,\cdots,\xi_{2(\mu+2)})$ still have the same meaning as defined in \eqref{fc1}.
We first consider the subbundle $\mathscr O(m_i)^{\oplus (\mu+1)}$ which has $(\xi_3,\xi_5,\xi_7,\cdots,\xi_{2\mu+3})$ as fiber coordinates.
Then by \eqref{hc2}, \eqref{hc3} and \eqref{fc1}, we obtain that 
 the image $\tilde\Phi(Z)$ is contained in a subvariety defined by
\begin{align}
\frac{u_2}{u_1}\xi_3=\xi_5,\,\,\frac{u_2}{u_1}\xi_5=\xi_7,\,\cdots, \,
\frac{u_2}{u_1}\xi_{2\mu+1}=\xi_{2\mu+3}.
\end{align}
These equations respectively define   a subbundle of $\mathscr O(m_i)^{\oplus \mu+1}$ of rank $\mu$, and as their intersection they define a line subbundle $\mathscr L$ of $\mathscr O(m_i)^{\oplus \mu+1}$.
By taking a meromorphic section of this line bundle concretely and counting the number of its zeros and poles, we can readily deduce that $\mathscr L\simeq\mathscr O(m_j)$, and that $\xi_3$ can be used as a fiber coordinate on $\mathscr L$.
By analogous computations for the subbundle $\mathscr O(m_i)^{\oplus (\mu+1)}$ having $(\xi_4,\xi_6,\xi_8,\cdots,\xi_{2\mu+4})$ as fiber coordinates, we obtain a line sub-bundle $\mathscr L'$ which is isomorphic to $\mathscr O(m_j)$ having $\xi_4$ as a fiber coordinate.
This way we deduce that $\tilde\Phi(Z)$ is contained in $\mathbf P(\mathscr O(m_i)^{\oplus 2}\oplus\mathscr L\oplus\mathscr L'\oplus\mathscr O)\simeq\mathbf P(\mathscr O(m_i)^{\oplus 2}\oplus\mathscr O(m_j)^{\oplus 2}\oplus\mathscr O)$.
Moreover, by \eqref{mt1} and \eqref{mt2} with $\alpha=2$, $\tilde\Phi(Z)$ satisfies the equations \eqref{de2}, as desired.
\proofend

\vspace{2mm}
In the sequel  the algebraic subvariety in $\mathbf P(\mathscr O(m_i)^{\oplus 2}\oplus\mathscr O(m_j)^{\oplus 2}\oplus\mathscr O)$ defined by \eqref{de2} is denoted by $X$, 
and let $p:X\to\Lambda_{m_i}$ be the projection.
$X$ is clearly a 3-dimensional irreducible subvariety and $T^2$ is acting preserving each fibers of $p$.
Hence $X$ is fibered  over $\Lambda_{m_i}$ by toric surfaces.
Recalling that  $\Lambda_{m_i}$ can be thought as a parameter space of the pencil $|F|^{T^2}$, it is tempting to expect that the last toric surfaces are bimeromorphic images of the members of the pencil (cf.\,\eqref{cd1}).
But the meromorphic map $\tilde\Phi:Z\to X$ is not necessarily bimeromorphic in general.
More precisely, $\tilde\Phi$ can be a generically $d:1$ map with $d>1$.
In order to determine the degree $d$, let $S\in |F|^{T^2}$ be any smooth real member, and for any $1\le \alpha\le k$, let $f_{\alpha}$ be a fiber of the holomorphic quotient map $\pi_{\alpha}:S\to \mathbf{CP}^1$ by the $G_{\alpha}^{\mathbf C}$-action which is the complexification of the $G_{\alpha}$-action on $S$.
Then we have  the following.

\begin{lemma}\label{lemma-degree}
Letting $d$ be the degree of the meromorphic map $\tilde\Phi:Z\to X$ as above,
$d=(f_i\cdot f_j)_S$ holds, where $(f_i\cdot f_j)_S$ means the intersection number in $S$.
\end{lemma}

\noindent
Since $(f_i\cdot f_j)_S>0$ always holds obviously, this in particular implies that $\tilde\Phi$ is surjective; namely $X=\tilde\Phi(Z)$ holds.
Further, the lemma also implies that $\tilde\Phi$ is bimeromorphic iff $(f_i\cdot f_j)_S=1$.
Hence $X$ is a projective model of $Z$ if $d=1$.
Note that if $j=i+1$, $(f_i\cdot f_j)_S=1$ holds. 
Note also that once $i$ and $j$ are given, it is easy to compute $(f_i\cdot f_j)_S$.

\vspace{1mm}
Proof of Lemma \ref{lemma-degree}. By the diagram \eqref{cd1}, the meromorphic map  $\tilde\Phi$ maps any two different members of $|F|^{T^2}$ to different fibers of $p$. 
Therefore, the degree of $\tilde\Phi$ is equal to the degree of the restriction $\tilde\Phi|_S$.
The last degree is of course equal to that of $\Phi|_S$.
The quotient map $\pi_{\alpha}$ has clearly two $T^2$-invariant fibers, which are reducible if $n>0$. We choose the fiber $f_{\alpha}$ as one of these $T^2$-invariant fibers, so that the conjugate curve $\ol f_{\alpha}$ is another $T^2$-invariant fiber.
Then since the divisor $Y_{\alpha}\in |m_{\alpha}F|$ (in Prop.\,\ref{prop-mt1}) satisfies $Y_{\alpha}|_S=m_{\alpha}C-f_{\alpha}+\ol f_{\alpha}$ as in \cite[Prop.\,2.5 (ii)]{Hon08-1}, $\Phi|_S$ is nothing but the meromorphic map associated to the 4-dimensional linear system generated by the following 5 effective curves
\begin{align}\label{5gen}
m_iC,\, m_iC-f_i+\ol f_i,\,m_iC+f_i-\ol f_i,\,
m_iC-f_j+\ol f_j, \,m_iC+f_j-\ol f_j.
\end{align}
Let $x_0,x_1,x_2,x_3,x_4\in H^0(m_iK_S^{-1})$ be sections defining the 5 curves \eqref{5gen} respectively.
We write $G=S\backslash C$, which is isomorphic to $\mathbf C^*\times\mathbf C^*$.
Then the 5 sections do not vanish on $G$.
Hence the ratio $x_a/x_0$ ($1\le a\le 4$) are non-constant and non-vanishing holomorphic functions on $G$.
By the definition of the meromorphic map associated to a linear system in general, on $G$, $\Phi|_S=(x_1/x_0,x_2/x_0,x_3/x_0,x_4/x_0)$ holds.
Now since $(x_1)-(x_0)=\ol f_i-f_i$ holds, the holomorphic function $x_1/x_0$ can be written in the form $\pi_i^*z$ (the pull-back of $z$ by $\pi_i$), where $z$ is an affine coordinate on the target space of $\pi_i$ centered at $\pi_i(\ol f_i)$. 
Similarly, since $(x_2)-(x_0)= f_i-\ol f_i$ holds, $x_2/x_0=\pi_i^*z^{-1}$ holds.
Parallel argument shows  $x_3/x_0=\pi_j^*w$ and $x_4/x_0=\pi_j^*w^{-1}$, where $w$ is an affine coordinate on the target space of $\pi_j$ centered at $\pi_j(\ol f_j)$.
Thus we obtain $\Phi|_S=(\pi_i^*z,\pi_i^*z^{-1},\pi_j^*w,\pi_j^*w^{-1})$.
This shows that an inverse image of $\Phi|_S$ of a point (in the image of $\Phi|_S$) is exactly the intersection of fibers of $\pi_i$ and $\pi_j$.
Thus the degree of $\Phi|_S$ is equal to $(f_i\cdot f_j)_S$, as desired.
\proofend


\vspace{2mm}
We reached the main result in this paper.

\begin{thm}\label{thm-1}
If $d=(f_i\cdot f_j)_S$ is one, 
the bimeromorphic map $\tilde\Phi:Z\to X$ is biholomorphic on the Zariski open subset $U\subset Z$ in \eqref{U}.
\end{thm} 

\noindent
Proof.
Since $\tilde\Phi$ is the composition of the map $\Phi:Z\to\mathbf P^{\vee}(W_{m_i}+V_{\mu}\cdot W_{m_j})$ and the inverse of the blow-up $\mathbf P(\mathscr O(1)^{\oplus (2\mu+4)}\oplus\mathscr O)\to\mathbf P(W_{m_i}+V_{\mu}\cdot W_{m_j})$ with center $\mathbf P_{\infty}$, it suffices to show that 
(i) $\Phi|_U:U\to\Phi(U)$ is biholomorphic, and 
(ii) $\Phi(U)\cap\mathbf P_{\infty}=\emptyset$.

First we show (i).
 We recall that $\Phi$ is $T^2$-equivariant and hence also $(\mathbf C^*)^2$-equivariant.
 By the diagram \eqref{cd1}, $\Phi$ can never maps two different 2-dimensional $(\mathbf C^*)^2$-orbits to one and the same $(\mathbf C^*)^2$-orbit.
Therefore,  since $\Phi$ is generically one to one by Lemma \ref{lemma-degree} under the assumption of the theorem, it suffices to show that $\Phi$ does not contract any 2-dimensional $(\mathbf C^*)^2$-orbits to  lower-dimensional orbits.
In the sequel we carefully verify this by distinguishing in which member of $|F|^{T^2}$ the 2-dimensional orbit is contained.

Suppose that a 2-dimensional $(\mathbf C^*)^2$-orbit $G$ is contained in an irreducible member $S\in |F|^{T^2}$. This means $S\backslash G=C$.
Then since $Y_l\cap S\subset C$ and $\ol Y_l\cap S\subset C$ for any $1\le l\le k$, the sections $y_i$, $\ol y_i$, $y_j$ and $\ol y_j$ do not have zero on $G$.
Therefore,  $g(u_1,u_2)y_j$ and $g(u_1,u_2)\ol y_j$ do not have zeros on $G$ for any monomial $g$ since $S\cap S_1^{\pm}\subset C$ and $S\cap S_2^{\pm}\subset C$.
Recalling that $\Phi$ is given by arranging the basis of $W_{m_i}+ V_{\mu}\cdot W_{m_j}$,
these imply that $\Phi(G)$ is 2-dimensional.

In the following suppose that $G\subset S_{\alpha}^+$ or $G\subset S_{\alpha}^-$ for some $1\le \alpha\le k$.

If both $Y_i$ and $Y_j$ do not contain $S_{\alpha}^+$ and $S_{\alpha}^-$, the 4 sections $y_i,\ol y_i,y_j$ and $\ol y_j$ do not have zeros on $G$. 
Hence at least one of $u_1^{m_j}y_j$ and $u_2^{m_j}y_j$ does not have zeros on $G$.
By the same reason as above, these imply that $\Phi(G)$ is 2-dimensional.

Suppose $G\subset S_{\alpha}^+\subset Y_i$.
Then $S_{\alpha}^+\not\subset \ol Y_i$ by Prop.\,\ref{prop-mt1} (iv).
Hence $\ol y_i$ does not have zeros on $G$.
Moreover, again by Prop.\,\ref{prop-mt1} (iv), at least one of $y_j$ and $\ol y_j$ does not have zeros on $G$.
If the former holds, at least one of  $u_1^{m_j}y_j$ and $u_2^{m_j}y_j$ does not have zeros on $G$.
With the above non-vanishing property for $\ol y_i$, this again implies that $\Phi(G)$ is 2-dimensional.
If the latter holds, at least one of $u_1^{m_j}\ol y_j$ and $u_2^{m_j}\ol y_j$ does not have zeros on $G$.
Again these imply that $\Phi(G)$ is 2-dimensional.
The case $G\subset S_{\alpha}^-\subset Y_i$ can be verified analogously.

Next suppose $G\subset S_{\alpha}^+\subset Y_j$.
Then by the same reason as above,  $\ol y_j$ does not have zeros on $G$.
Hence at least one of $u_1^{m_j}\ol y_j$ and $u_2^{m_j}\ol y_j$ has no zeros on $G$.
On the other hand, at least one of $y_i$ and $\ol y_i$ has no zeros on $G$.
Thus we again obtain that $\Phi(G)$ is 2-dimensional.
The case $G\subset S_{\alpha}^-\subset Y_j$ can be verified analogously.
This way we have covered all cases, and conclude that $\Phi$ is biholomorphic on any 2-dimensional orbits. 
%
Thus we obtain the claim (i).
(The above argument shows that $\Phi$ does not contract any two-dimensional $(\mathbf C^*)^2$-orbits even when $\Phi$ fails to be bimeromorphic.)

Finally we show the claim (ii).
As in \eqref{center1} the center $\mathbf P_{\infty}$ is explicitly given by $\{\zeta_{N-m_i}=\cdots=\zeta_N =0\}$. 
By \eqref{hc3}, we have 
\begin{align}
\Phi^{-1}(\mathbf P_{\infty})&=\{z\in Z\set u_1(z)^{m_i}=u_1(z)^{m_i-1}u_2(z)=\cdots=u_2(z)^{m_i}=0\}\\
&=\{z\in Z\set u_1(z)=u_2(z)=0\}.\label{c2}
\end{align}
But since $C={\rm{Bs}}\,|F|^{T^2}$ and $H^0(F)^{T^2}=\langle y_1,y_2\rangle$, \eqref{c2} is precisely $C$.
Therefore, because $U\cap C=\emptyset$, we have $\Phi(U)\cap\mathbf P_{\infty}=\emptyset$, as desired.
\proofend

\vspace{2mm}
Thus provided $(f_i\cdot f_j)=1$, it has turned out  that the projective algebraic variety $X$ contains  the twistor space $Z^{\circ}$ of the explicit Joyce metric on 
$\mathscr H^2\times T^2$.
Of course, $Z^{\circ}$  is a non-Zariski dense open subset of $X$.
Furthermore, by the above proof of the claim (ii), it is obvious that 
\begin{align}\label{affine}
\tilde{\Phi}(U)\subset \{(\xi_1,\xi_2,\xi_3,\xi_4;\lambda)\in\mathscr O(m_i)^{\oplus 2}\oplus\mathscr O(m_j)^{\oplus 2}\set \hspace{1mm}{\text{subject to the equations}}\hspace{2mm} \eqref{de2}\}
\end{align} 
where the rank-4 vector bundle in \eqref{affine} is considered as a Zariski open subset of the $\mathbf P^4$-bundle $\mathbf P(\mathscr O(m_i)^{\oplus 2}\oplus\mathscr O(m_j)^{\oplus 2}\oplus\mathscr O)\to\Lambda_{m_i}$ as before.
Thus the twistor space $Z^{\circ}$ is realized algebraically in the rank-4 vector bundle over $\mathbf{CP}^1$.
In particular, $Z^{\circ}$ can be holomorphically embedded in a projective space.

The structure of fibers of the projection $p:X\to\Lambda_{m_i}$ is readily determined from the defining equations.
All fibers are  toric surfaces and are complete intersections of two hyperquadrics in $\mathbf{CP}^4$. 
If both of the right-hand side of \eqref{de2} are non-zero for $\lambda\in\Lambda_{m_i}$, the fiber $p^{-1}(\lambda)$ is irreducible (toric) surface which has 4 ordinary double points.
These are birational images of smooth members of the pencil $|F|^{T^2}$.
If just one of the right-hand side vanishes,
$p^{-1}(\lambda)$ splits into two quadratic cones of full rank in two hyperplanes $\mathbf{CP}^3$.
Both of these cones are birational images of a reducible member of $|F|^{T^2}$.
If both of the right-hand side vanish, $p^{-1}(\lambda)$ splits into four linear planes.
Precisely two of the planes are birational image of a reducible member of $|F|^{T^2}$.
We note that if $Y_i$ or $Y_j$ has non-reduced components, the projective model $X$ has singularities along a curve contained in a reducible fiber of $p$.
This kind of singularity does not occur in \cite{Hon06-1}, and makes  it difficult to construct the biholomorphic models of the twistor spaces in general.





\small
\vspace{13mm}
\hspace{7.5cm}
$\begin{array}{l}
\mbox{Department of Mathematics}\\
\mbox{Graduate School of Science and Engineering}\\
\mbox{Tokyo Institute of Technology}\\
\mbox{2-12-1, O-okayama, Meguro, 152-8551, JAPAN}\\
\mbox{{\tt {honda@math.titech.ac.jp}}}
\end{array}$

\end{document}